\documentclass[12pt]{amsart}

\usepackage{amssymb}
\usepackage{amsfonts}
\usepackage{amsmath}
\usepackage{amssymb}
\usepackage[mathcal]{eucal}
\usepackage{amscd}
\usepackage{verbatim}

\newtheorem{thm}{Theorem}[section]
\newtheorem{cor}[thm]{Corollary}
\newtheorem{lem}[thm]{Lemma}
\newtheorem{prop}[thm]{Proposition}

\newtheorem{question}[thm]{Question}
\theoremstyle{definition}
\newtheorem{defn}[thm]{Definition}

\theoremstyle{remark}
\newtheorem{remark}[thm]{Remark}

\numberwithin{equation}{section}

\newcommand{\Acal}{\mathcal{A}}

\newcommand{\Dcal}{\mathcal{D}}

\newcommand{\Fcal}{\mathcal{F}}

\newcommand{\Pcal}{\mathcal{P}}

\newcommand{\Scal}{\mathcal{S}}

\newcommand{\Ucal}{\mathcal{U}}

\newcommand{\Xcal}{\mathcal{X}}
\newcommand{\Ycal}{\mathcal{Y}}

\newcommand{\Q}{\mathbb{Q}}
\newcommand{\R}{\mathbb{R}}
\newcommand{\C}{\mathbb{C}}
\newcommand{\N}{\mathbb{N}}

\newcommand{\ip}[2]{\langle#1,#2\rangle}

\newcommand{\al}{\alpha}
\newcommand{\Ga}{\Gamma}
\newcommand{\ga}{\gamma}
\newcommand{\del}{\delta}
\newcommand{\Del}{\Delta}
\newcommand{\ep}{\varepsilon}
\newcommand{\sig}{\sigma}

\newcommand{\om}{\omega}
\newcommand{\Om}{\Omega}

\newcommand{\imp}{\ \Rightarrow\ }

\newcommand{\tri}{\bigtriangleup}
\newcommand{\rest}{\upharpoonright}

\newcommand{\inte}{\rm{int\,}}

\newcommand{\Aut}{\rm{Aut\,}}

\newcommand{\Homeo}{\rm{Homeo\,}}
\newcommand{\MALG}{\rm{MALG\,}}
\newcommand{\Lip}{{\rm{Lip\,}}}

\newcommand{\sgn}{{\rm{sgn\,}}}

\newcommand{\diam}{\rm{diam\,}}
\newcommand{\supp}{\rm{supp\,}}

\newcommand{\ess}{\rm{ess}}

\newcommand{\tr}{\rm{tr\,}}

\newcommand{\dist}{{\rm{dist}}}

\newcommand{\VarLip}{{\rm{VarLip\,}}}

\renewcommand{\(}{\bigl(}
\renewcommand{\)}{\bigr)\vphantom{)}}

\newcommand{\cP}[2]{\mathbb{P}\,\(\,#1\,\big|\,#2\,\)\,}

\newtheorem{rmk}[thm]{Remark}

\begin{document}

\title
[Automorphisms of Gaussian measure]
{The automorphism group of the Gaussian measure
cannot act pointwise}

\author{E. Glasner, B. Tsirelson and B. Weiss}

\address{Department of Mathematics\\
     Tel Aviv University\\
         Ramat Aviv\\
         Israel}
\email{glasner@math.tau.ac.il}

\address{Department of Mathematics\\
     Tel Aviv University\\
         Ramat Aviv\\
         Israel}
\email{tsirel@post.tau.ac.il \hfill
 http://www.tau.ac.il/\textasciitilde tsirel/}

\address {Institute of Mathematics\\
 Hebrew University of Jerusalem\\
Jerusalem\\
 Israel}
\email{weiss@math.huji.ac.il}

\begin{abstract}
Classical ergodic theory deals with measure (or measure class)
preserving actions of locally compact groups on Lebesgue spaces.  An
important tool in this setting is a theorem of Mackey which provides
spatial models for Boolean $G$-actions. We show that in full
generality this theorem does not hold for actions of Polish groups. In
particular there is no Borel model for the Polish automorphism group
of a Gaussian measure. In fact, we show that this group as well as
many other Polish groups do not admit any nontrivial Borel measure
preserving actions.
\end{abstract}

\maketitle

\section*{Introduction}

Our motivation is threefold: invariant measures; Borel liftings;
Gaussian measures.

\begin{description}
\item[Invariant measures]
By a famous theorem of A. Weil if a Polish group $ G $ admits a
$\sig$-finite invariant measure then $ G $ is locally compact
(see Appendix B).

Nonetheless, even if $G$ is not locally compact, a homogeneous space
of $G$ might even admit a finite invariant measure.  For example, the
group $G$ of $\mu$-preserving homeomorphisms of the Cantor set
$\Omega$, will act transitively on $\Omega$ for a suitable choice of
$\mu$ (see for example \cite{GW}).  We show that this never happens
for some classes of Polish groups $ G $ (for instance, the full
unitary group of a separable Hilbert space), except for the trivial
case: a measure concentrated on fixed points.

\item[Borel liftings]
Let $ G $ be a closed subgroup of the Polish group of all invertible
measure preserving transformations of (say) $ [0,1] $ with Lebesgue
measure. An element $ g \in G $ is an equivalence class of maps $
[0,1] \to [0,1] $ rather than a single map; thus, $ g(x) $ is defined
only almost everywhere. Can we define $ g(x) $ everywhere?  More
exactly: can we lift the $ \bmod\,0 $ action to a Borel action? We
give a general criterion for lifting, and a negative answer for some
classes of groups including the Gaussian case.

\item[Gaussian measures]
Every Euclidean space carries its standard Gaussian measure. However,
a separable Hilbert space $ H $ (over $ \R $) does not. The standard
Gaussian process over $ H $ is a linear isometry between $ H $ and the
subspace spanned by a sequence of i.i.d.\ $ N(0,1) $ random variables
in $ L_2 $ over a probability space. Can we implement each point of
(some version of) the probability space as a function on $ H $ or
another superstructure over $H$?  Some well-known such constructions
are `isotropic', that is, invariant under the full orthogonal group of
$ H $. Others give a standard probability space.  We show that these
two desirable properties exclude each other.  Two proofs are given,
one via `invariant measures', the other via `Borel liftings'.

\end{description}

Having thus stated our goals in outline let us be more precise.
Traditionally, ergodic theory is treated within the context of locally
compact groups acting on standard Lebesgue probability spaces. However
it is often the case that one has to deal with near-actions (see
definition below) or merely with an action of the group on a measure
algebra (i.e. the Borel algebra modulo sets of measure zero) and it is
then desirable to find a standard Lebesgue model, or even better, a
Polish (= complete metric, second countable space) or a compact model
where the group acts continuously.

Recall that a \emph{Borel action} of $ G $ on a Borel space
$(X,\Xcal)$ is a Borel map $ G \times X \to X $ (denote it just $
(g,x) \mapsto gx $) satisfying the two conditions, $ ex = x $ and $
g(hx) = (gh)x $ for all $ g,h \in G $ and all $ x \in X $. Such an
object is called also a \emph{Borel $ G $-space}.

\begin{defn}\label{def-near}
(See Zimmer \cite[Def.~3.1]{Z})
Let $G$ be a Polish group and $(X,\Xcal,\mu)$ a standard Borel space
with a probability measure $\mu$.  By a \emph{near-action\/} of $G$ on
$(X,\Xcal,\mu)$ we mean a Borel map $G\times X\to X, (g,x)\mapsto gx$
with the following properties:
\begin{enumerate}
\item[(i)]
With $e$ the identity element of $G$, $ex=x$ for almost every $x$.
\item[(ii)]
For each pair $g,h \in G, g(h x)=(gh) x$ for almost every $x$
(where the set of points $x\in X$ of measure one where this
equality holds may depend on the pair $g, h$).
\item[(iii)]
Each $g\in G$ preserves the measure $\mu$.
\end{enumerate}
\end{defn}

Let $ \Aut(X)= \Aut(X,\Xcal,\mu) $ be the Polish group of all
equivalence classes of invertible measure preserving transformations
$X \to X$,
with the neighborhood basis at the identity formed by
sets of the form
$$
N(A,\ep)=\{T\in {\Aut}(X): \mu(A\triangle TA)<\ep\},
$$
for $A\in \Xcal$ and $\ep>0$.
What we would like to show next is that the following three notions
are equivalent.

\begin{enumerate}
\item[(I)]
A near-action of $ G $ on $(X,\Xcal,\mu)$.
\item[(II)]
A continuous homomorphism from $ G $ to $ \Aut(X) $.
\item[(III)]
A Boolean action of $ G $ on $(X,\Xcal,\mu)$, that is, a continuous
homomorphism from $ G $ to the automorphism group of the associated
measure algebra.
\end{enumerate}

Given a near action of $G$, it is easy to check that the
natural mapping from $G$ to $\Aut(X)$ defines a measurable mapping.
That it is a homomorphism follows from the defining property of being
a near action, and since, as is well known, measurable homomorphisms
of Polish groups are continuous, we get (II) from (I). To go in the
other direction, we must construct from a continuous homomorphism of
$G$ into $\Aut(X) $, a near action of $G$ on $(X,\Xcal,\mu)$.

For this we need to define a section on equivalence classes of Borel
measurable functions where the equivalence relation is that of
equality $\mu$ a.e. Let $ (X,\Xcal) $ and $ (Y,\Ycal) $ be standard
Borel spaces and $ \mu $ a probability measure on $(X,\Xcal) $. Then
the set $ L_0 (X,Y) = L_0 ( (X,\Xcal,\mu), (Y,\Ycal) ) $ of all
equivalence classes ($\bmod\,0$ with respect to $\mu$) of Borel (or
just $\mu$-measurable) maps $ X \to Y $ is also a standard Borel
space; its $\sig$-algebra is generated by functions $ f \mapsto \mu(A
\cap f^{-1}(B)) $ for $ A \in \Xcal $, $ B \in \Ycal $.

There exists a (highly non-unique) Borel map $ V : L_0(X,Y) \times
X \to Y $ such that for every $ f \in L_0(X,Y) $ the function $ x
\mapsto V(f,x) $ belongs to the equivalence class $ f $. For
example, assuming $ X=Y=(0,1) $ (with the usual Borel
$\sig$-algebra and Lebesgue measure $\mu$), we may take
\[
V(f,x) = \limsup_{\ep\to0} \frac1{2\ep} \int_{x-\ep}^{x+\ep}
f(x_1) \, dx_1 \, .
\]

The Polish group $ \Aut(X)$ is a Borel subspace of $ L_0(X,X) $.
Returning now to our situation, let $\phi$ denote a continuous
homomorphism of $G$ into $\Aut(X)$.  Composing $\phi$ with the
restriction of $V$ above to $ \Aut(X)\times X $ gives us a Borel
mapping from $G\times X$ to $X$, and one checks easily that the
properties for being a near action are satisfied. Thus (II) implies
(I).

Finally, $ \Aut(X) $ may be thought of as the automorphism group of
the measure algebra $ \MALG(X,\Xcal,\mu) = (\bar\Xcal,\bar\mu)$, where
$\bar\Xcal$ is $\Xcal$ modulo nullsets, and $\bar\mu$ the
corresponding measure; automorphisms of the measure algebra must
preserve Boolean operations and $\bar\mu$. Thus (II) and (III) are
equivalent.

This completes the discussion of the equivalence of the various
notions of a near-action.

In contrast to near-actions, we define the notion of spatial action.

\begin{defn}\label{def-spatial}
Let $G$ be a Polish group. By a \emph{spatial $G$-action } we mean a
Borel action of $G$ on a standard Lebesgue space $(X,\Xcal,\mu)$
such that each $g\in G$ preserves the measure $\mu$.  We say that two
spatial actions are isomorphic, if there exists a measure preserving
one to one map between two $G$-invariant subsets of full measure in
the corresponding spaces which intertwines the $G$-actions (the same
two sets for all $g\in G$).
\end{defn}

Every spatial action is also a near-action. In that case the spatial
action will be called a \emph{spatial model} of the near-action (or
the corresponding Boolean action). The question is, when does a given
near-action admit a spatial model.

Rohlin, in \cite{Ro}, when discussing $\R$-flows, distinguishes
between these two notions calling our near-actions continuous flows
and the spatial actions measurable flows. He notes there that the
theories of these two notions are not equivalent.
Indeed already J. von Neumann in his foundational work \cite{vN}
was aware of these distinctions and in footnote 13 writes that he
hopes to provide a proof that every near-flow has a continuous
spatial model.

We recall that a Polish $G$-space is a Polish space $ X $ together
with a continuous action $ G \times X \to X $ of a Polish group $ G
$. Such an action will be called a Polish action. If in addition $ X $
is compact then it is a compact Polish $ G $-space.

Every Polish action is also a Borel action. In that case the Polish
action will be called a \emph{Polish model} of the Borel action.

We have the following classical theorems, due to Mackey, Varadarajan
and Ramsay (\cite[Th.~3.2]{Va0}, \cite{Ma},
\cite[Th.~3.3]{Ra} and \cite{Va}; see also \cite{Ve}
and \cite{D}).

\begin{thm}\label{lc-model}
Let $G$ be a locally compact second countable topological group.

\textup{(a)}
Every near-action \textup{(}or Boolean action\textup{)} of $ G $
admits a spatial model.

\textup{(b)}
Every spatial action of $ G $ admits a Polish model.
\end{thm}

A powerful generalization to Polish groups of Theorem
\ref{lc-model}(b), given in \cite[Th.~5.2.1]{BK}, is crucial for our
work.

\begin{thm}[Becker and Kechris]\label{BK-thm}
\mbox{ }

\textup{(a)}
Every Borel action of a Polish group admits a Polish model.

\textup{(b)}
Every Borel $G$-space is embedded \textup{(}as a $G$-invariant Borel
subset\textup{)} into a compact Polish $G$-space.
\end{thm}

Item (b) above follows from \cite[Th.~2.6.6]{BK} (which in turn
utilizes a theorem of Beer \cite{B}).

In the present work we show that a full generalization to Polish
groups of the first part of Theorem \ref{lc-model} is not possible.
Many near-actions (or Boolean actions) of Polish groups admit no
spatial models.

In Section 1 we recall the definition of L\'{e}vy groups,
a class of groups which includes $\Aut(X,\Xcal,\mu)$ and
$U(H)$, the unitary group on an infinite dimensional Hilbert space.
We show that these groups admit no non-trivial spatial actions,
and discuss some further examples. In Section 2 we take up the
general question of finding criteria for a near-action to admit
a spatial model, or more generally a spatial factor.
We give a necessary and sufficient condition for this in terms of
$G$-continuous functions. These are those functions
$f \in L^\infty(\mu)$ with the property that $f \circ
g_n$ converges to $f$ in $L^\infty(\mu)$ norm whenever $g_n \to e$,
the unit element in $G$.

In order to apply the criterion of Section 2 we introduce in Section 3
whirly actions which may be defined as near-actions of
$G$ on $(X,\Xcal,\mu)$ such that for all
sets $A,B\in\Xcal$ of positive measure, every neighborhood of $e$
contains $g$ with $ \mu (A\cap gB) > 0 $.
We show that whirly actions have no non-constant $G$-continuous
functions and verify some easy examples of whirly actions.

In Section 4 we return to the natural near-action
of the Polish orthogonal group $G=O(H)$, when $H$ is the ``first
chaos'' Gaussian Hilbert space and show that this near-action is
whirly, thus giving the promised second proof of the non-existence
of a spatial model for this action.

Appendix A collects elementary proofs of the
L\'{e}vy property for many well known L\'{e}vy families.
In Appendix B we provide a proof of the fact,
mentioned at the beginning of the introduction,
that a Polish group which admits a
Borel $\sig$-finite invariant measure is locally compact.

We thank A. Vershik, V. Milman and A. Kechris for
instructive remarks concerning the historical background
of this work. We thank V. Pestov for a careful reading of
a draft of this paper and for several useful remarks.

\section{L\'evy groups admit no spatial actions}

The phenomenon of concentration of mass was first
considered by E. Borel in his ``law of large numbers",
where it is manifested in the family
of ``discrete cubes" $\{0,1\}^n$ equipped with Hamming distance
and counting measure. In the 30's P. L\'evy studied
the concentration phenomenon for the family of
Euclidean spheres. Then about 1970 V. Milman in
\cite{Mi1} revitalized the area when he discovered a new proof
of Dvoretzky's theorem using the concentration
phenomenon on spheres. In the work \cite{GM} M. Gromov and V. Milman
considered applications of the concentration phenomenon in
topological dynamics. In particular the notion of a
L\'evy family is introduced in \cite{GM}. See \cite{Mi4}
and \cite{P1} for further details on the history of
this subject.

Let $(X_n,d_n,\mu_n)$, $n=1,2,3\dots$ be a family of metric spaces
with probability measures $\mu_n$. Call such a family a \emph{L\'evy
family} if the following condition is satisfied.  If $A_n\subset X_n$
is a sequence of subsets such that $\liminf \mu_n(A_n)>0$ then for any
$\ep >0$, $\lim \mu(B_\ep (A_n))=1$, where $B_\ep(A)$ is the $\ep$
neighborhood of $A$.

A Polish group $G$ is a \emph{L\'evy group\/} if there exits a family
of compact subgroups $K_n\subset K_{n+1}$ such that the group
$F=\cup_{n\in \N} K_n$ is dense in $G$ and the corresponding family
$(K_n,d,m_n)$ is a L\'evy family; here $m_n$ is the normalized Haar
measure on $K_n$, and $ d $ is a right-invariant compatible metric on
$ G $ (the choice of $ d $ does not matter). Using left-invariant
metrics instead, we get an equivalent definition (just apply the map $
g \mapsto g^{-1} $).

Here is a list of some Polish groups well-known to be L\'evy groups.
Refer to Milman \cite{Mi2}, \cite{Mi3},
Gromov and Milman \cite{GM}, Glasner \cite{G},
Pestov \cite{P1} and Giordano and Pestov \cite{GP};
also see Appendix A for more details.

\begin{itemize}
\item
The full unitary group $ U(H) $ (of a separable Hilbert space $ H $);
and the full orthogonal group $ O(H) $ (of a separable Hilbert space
$ H $ over $ \R $),
where both groups are equipped with the
strong operator topology.
\item
The dense subgroup of $ U(H) $ consisting of
all unitary (or orthogonal) $ U $ such that $
\tr ((1-U)^*(1-U)) < \infty $.
\item
The group $ \Aut(X) $ (mentioned in the Introduction).
\item
The commutative (moreover, monothetic) group $ L_0 \( [0,1], S^1 \) $
of all (equivalence classes of) measurable functions $ [0,1] \to S^1
$, where $ S^1 = \{ z \in \C : |z| = 1 \} $.
\end{itemize}

\begin{thm}\label{thm-levy}
Every spatial action of a L\'{e}vy group is trivial; i.e. the set of
fixed points is of full measure.
\end{thm}

\begin{proof}
By Theorem \ref{BK-thm}(b), every Borel $G$-space is embedded into
a compact Polish $G$-space. Therefore it suffices to prove the theorem
for a continuous action of $ G $ on a metrizable compact space $ X $
and a $G$-invariant Borel probability measure $ \mu $. We will see that
$ G $ acts trivially on $ \supp\mu $ (the support of the measure).

The action is a continuous homomorphism from $ G $ to the Polish group
$ {\Homeo}(X) $ of all homeomorphisms of $ X $ (as noted in
\cite[p.~427]{P1}). We equip $ X $ with a compatible metric $ \rho $,
and $ {\Homeo}(X) $ with the compatible right-invariant metric $ (f,g)
\mapsto \max_{x\in X} \rho(f(x),g(x)) $. Now the homomorphism is
uniformly continuous, provided that $ G $ is also equipped with a
right-invariant metric (which will be assumed).

The family $ ( g \mapsto g \cdot x )_{x\in X} $ of maps $ G \to X $ is
equicontinuous. By \cite[2.1]{GM}, it sends the L\'evy family $
(m_n) $ of measures on $ G $ to a L\'evy family $ (m_n \cdot
x) $ of measures on $ X $, uniformly in $ x \in X $. In other words:
for all $ x_1, x_2, \dots \in X $ the family $ (m_n \cdot x_n) $
is L\'evy.

By \cite[2.4]{GM} the family of measures is degenerate in the sense
that
\[
\min_{y\in X} \int_X \rho(\cdot,y) \, d(m_n\cdot x) = \min_{y\in
X} \int_{K_n} \rho(g\cdot x,y) \, dm_n(g) \to 0
\]
for $ n \to \infty $, uniformly in $ x \in X $. The proof is
simple. Assuming the contrary and using compactness, we choose $ x_k
\in X $ and $ n_k \to \infty $ such that measures $ m_{n_k} \cdot x_k
$ converge (weakly) to some measure $ \nu $ on $ X $ satisfying $
\min_{y\in X} \int_X \rho(\cdot,y) \, d\nu > 0 $, which means that the
support of $ \nu $ contains at least two points. Every open set $ A
\subset X $ such that $ \nu(A) > 0 $ satisfies $ \nu(B_\ep(A)) = 1 $
for all $ \ep > 0 $; here $ B_\ep(A) $ is the closed $ \ep
$-neighborhood of $ A $. We get a contradiction by choosing $ A $ such
that some points of the support of $ \nu $ belong to $ A $ and some do
not belong to the closure of $ A $.

For each $ n $ we introduce the subspace
$ H_n \subset L_2(\mu) $ of all $ K_n $-invariant functions, and the
corresponding orthogonal projection $ Q_n $,
\[
Q_n f(x) = \int_X f \, d(m_n\cdot x) = \int_{K_n} f(g\cdot x) \,
dm_n(g) \, .
\]
For $ f \in C(X) \subset L_2(X) $ the degeneracy of measures gives
\[
\int | f(g\cdot x) - (Q_n f)(x) |^2 \, dm_n(g) \to 0
\]
for $ n \to \infty $, uniformly in $ x $.

We see this as follows. Denoting by $y_{n,x} $ the minimizer of
$ \int_X \rho(\cdot,y) \, d(m_n\cdot x) $
we have
\begin{multline*}
\int | f(g\cdot x) - (Q_n f)(x) |^2 \, dm_n(g) = \min_{a\in\R} \int |
 f(\cdot) - a |^2 \, d(m_n\cdot x) \\
\le \int | f(\cdot) - f(y_{n,x}) |^2 \, d(m_n\cdot x) \le
 \int_{B_\ep(y_{n,x})} + \int_{X\setminus B_\ep(y_{n,x})} \\
\le \max_{B_\ep(y_{n,x})} | f(\cdot) - f(y_{n,x}) |^2 + \Big( \max_X |
 f(\cdot) - f(y_{n,x}) |^2 \Big) \frac1\ep \int_X \rho(\cdot,y_{n,x})
 \, d(m_n\cdot x) \, .
\end{multline*}
Taking $ \limsup_{\ep\to0} \limsup_{n\to\infty} \sup_{x\in X} (\dots)
$ we get $ 0 $.

On the other hand, the
integral
\[
\int | f(g\cdot x) - (Q_n f)(x) |^2 \, d\mu(x)
\]
does not depend on $ g \in K_n $ and is equal to $ \| f - Q_n f \|^2
$. Therefore
\[
\| f - Q_n f \|^2 = \iint | f(g\cdot x) - (Q_n f)(x) |^2 \,
dm_n(g) \, d\mu(x) \to 0
\]
for $ n \to \infty $. However, the inclusion $ K_n \subset K_{n+1} $
implies $ H_n \supset H_{n+1} $ and $ \| f - Q_n f \| \le \| f -
Q_{n+1} f \| $. So, $ \| f - Q_n f \| = 0 $. It means that $ f $ is $
K_n $-invariant, that is, $ f(x) = f(g\cdot x) $ for all $ x \in
\supp\mu $ and all $ f \in C(X) $. Thus, $ g\cdot x = x $ for all such
$ x $ and all $ g \in \cup K_n $, therefore all $ g \in G $.
\end{proof}

\begin{question}
Can a L\'evy group admit a nontrivial nonsingular \textup{(}that is,
preserving a measure class\textup{)} Borel action\textup{?}
\end{question}

\begin{remark}
The basic idea in the proof of Theorem \ref{thm-levy} is derived from
Gromov and Milman \cite{GM} where they show that L\'evy groups have
the fixed point on compacta property. The question arises whether
every group with the fixed point property do not admit a nontrivial
spatial measure preserving Borel action. Now it was shown by Pestov
\cite{P} that the Polish group $G=\Aut(\Q,<)$ of order preserving
permutations of the rational numbers, equipped with the topology of
pointwise convergence (with respect to the \emph{discrete} topology on
$\Q$), has the fixed point on compacta property (or is extremely
amenable).  However it is easy to see that this group also acts
ergodically by homeomorphisms on the ``$\Q$-Bernoulli system''
$(\Om,\Fcal,\mu)$. Here $\Om=\{1,-1\}^{\Q}$,\ $\mu$ is the product
measure $\mu=(1/2,1/2)^{\Q}$ and $G$ acts on ``configurations''
$\om\in \Om$ by permuting the indices. We therefore conclude that some
Polish groups with the fixed point property can have nontrivial
spatial actions.
\end{remark}

\begin{remark}
Recall that a topological group $G$ is \emph{amenable}
if each compact $G$-space admits a $G$-invariant probability
measure. Using this definition of amenability and the fact
that every compact group is amenable it is easy to deduce that
every L\'evy group is amenable. Now if $G$ is a Polish L\'evy
group and $(X,G)$ is a compact $G$-space, then by amenability
of $G$ there is a $G$-invariant probability measure
$\mu$ on $X$. It can be shown that $(X,G)$ is
represented as an inverse limit of a directed system of metrizable
$G$-spaces $\{(X_\al,G)\}$. Let $\mu_\al$ be the
image of $\mu$ on $X_\al$, then apply theorem \ref{thm-levy}
to deduce that ${\supp}(\mu_\al)$ is a closed nonempty
collection of fixed points. It is now easy to conclude that
the support of $\mu$, ${\supp}(\mu)$ is a nonempty closed
subset of $X$ consisting of fixed point. Thus the Gromov-Milman
theorem that every L\'evy group has the fixed point on compacta
property follows from Theorem \ref{thm-levy}. Of course
this is a rather circumventive way of proving it.
\end{remark}

\begin{remark}
The following application of Theorem \ref{thm-levy} was
pointed out to us by V. Pestov. Some years ago he and
M. Cowling conjectured that
every invariant mean on the unitary group $ U(H) $ is
contained in the weak$^*$ closed convex hull of the
multiplicative invariant means. Now for any topological group
$G$ the above statement holds iff every $G$-invariant
measure on the greatest ambit $\Scal(G)$ of $G$
(i.e. the Gelfand space of the Banach algebra $BLUC(G)$
of bounded left uniformly continuous functions on $G$)
is supported on the set of fixed points. Again
using the fact that for a Polish L\'evy group $\Scal(G)$ is an
inverse limit of a directed system of metrizable
$G$-spaces, we deduce from Theorem \ref{thm-levy} that
\emph{every invariant mean on a Polish L\'evy group is
contained in the weak$^*$ closed convex hull of the
multiplicative invariant means}.
\end{remark}

\begin{remark}
Note that, for example, the group $S_\infty$ of permutations of $\N$
(with the topology of pointwise convergence) is a non locally compact
Polish subgroup of $\Ucal(H)$ which admits nontrivial measure
preserving spatial actions.
\end{remark}

\begin{remark}
There are well known examples of Polish groups $G$ which do not admit
any weakly continuous linear representations on a Banach space; see
e.g.  \cite{HC} and \cite{Ba}. In the latter Banaszcyk provides, for
every infinite dimensional normed space $E$, examples of the form
$G=E/K$ where $K\subset E$ is a discrete subgroup. It is easy to see
that any such group is moreover monothetic. Of course such ``strongly
exotic groups'' as they are called by Herer, Christensen and Banaszcyk
can not admit even a nontrivial near-action.  Moreover, every
nonsingular near-action (preserving a measure class rather than a
measure) leads, by a standard construction, to a unitary
representation. Thus these strongly exotic groups can not admit
nontrivial \emph{nonsingular} near-actions.
\end{remark}

By Theorem \ref{thm-levy}, a nontrivial near-action of a L\'evy group
cannot admit a spatial model.  An important example is the
automorphism group of an infinite-dimensional Gaussian measure. Up to
isomorphism, the relevant probability space is the product $
(\R^\infty,\ga^\infty) $ of countably many copies of $ (\R,\ga) $,
where $ \ga $ is the standard one-dimensional Gaussian measure (normal
distribution). The space (so-called \emph{first chaos}) of all
measurable linear functionals on $ (\R^\infty,\ga^\infty) $ is $ l_2
$. The action of the full orthogonal group $ O(l_2) $ on measurable
linear functionals is well-known to be induced by its near-action on $
(\R^\infty,\ga^\infty) $, which is what we mean by the automorphism
group of the Gaussian measure. (In this sense, $ O(l_2) $ is a closed
subgroup of $ \Aut(\R^\infty,\ga^\infty) $; see also Section 4.)

\begin{cor}\label{nogo1}
The near-action of the automorphism group of the Gaussian measure
admits no spatial model.
\end{cor}

\begin{remark}
In \cite[Theorem 4]{Ve} A. Vershik states:
There is no measurable realization of the group $U(H)$, that is,
there is no set of full measure that is invariant under all
$u\in U(H)$. This would yield Corollary 1.6, however the proof
given there appears to us to be incomplete.
\end{remark}

\begin{remark}
Another proof of Corollary \ref{nogo1} uses the L\'evy group $ L_0 \(
[0,1], S^1 \) $ rather than $ O(l_2) $. The latter group contains (an
isomorphic copy of) the former group as a closed subgroup, see
e.g. Lema\'nczyk, Parreau and Thouvenot \cite{LPT}. The near-action of
(the copy of) $ L_0 \( [0,1], S^1 \) $ on $ (\R^\infty,\ga^\infty) $
is nontrivial; by Theorem \ref{thm-levy} it cannot admit a spatial
model, which implies Corollary \ref{nogo1}.
\end{remark}

\begin{remark}
About the meaning of Corollary \ref{nogo1}.
Almost all points of $ (\R^\infty,\ga^\infty) $ do not belong to $ l_2
$ and therefore cannot be interpreted as continuous linear functionals
on $ l_2 $. One could hope for interpreting them as another
superstructures over $ l_2 $ (say, densely defined discontinuous
linear functionals) that form a Borel $ G $-space ($ G $ being the
symmetry group). Corollary \ref{nogo1} shows that it is
impossible. Maybe, Borel measurability could be weakened (say, to
universal measurability)? We do not know. Some related ideas can be
found in \cite[Prop.~E.2]{Du} and \cite[Example 1.27]{Ja}.
\end{remark}

\section{Which actions admit spatial models ?}

In this section we enhance our understanding of the lifting problem by
relating it to a notion of $G$-continuity of functions which is
reminiscent of the classical notion of a rigid action in ergodic
theory.

\begin{defn}
Having a near-action (or Boolean action) of $G$ on $(\Xcal,\mu)$ we
say that $f \in L^\infty(\mu)$ is \emph{$G$-continuous\/}, if $f \circ
g_n$ converges to $f$ in $L^\infty(\mu)$ norm whenever $g_n \to e$.
\end{defn}

The collection $\Acal(G)$ of all $G$-continuous functions is a
$G$-invariant closed subalgebra of $L^\infty(\mu)$.

\begin{thm}\label{2.2}
A near-action admits a spatial model if and only if there exists a
sequence of $G$-continuous functions that generates the $\sig$-algebra
(equivalently: separates points).
\end{thm}

\begin{proof}
Suppose first that we have a spatial model, that is, a Borel $G$-space
with an invariant measure. By Theorem \ref{BK-thm}(b) this Borel
$G$-space can be embedded into a compact Polish $G$-space $X$ (with an
invariant measure). The continuous functions on $X$ form a separable
Banach space and a dense sequence in $C(X)$ will provide a sequence of
$G$-continuous functions in $L^\infty(\mu)$ which separates points.

Conversely, suppose there exists a sequence $\{f_n:n\in \N\}\subset
L^\infty(\mu)$ of $G$-continuous functions that generates the
$\sigma$-algebra. Let $G_0\subset G$ be a countable dense subgroup of
$G$.  Let $A\subset L^\infty(\mu)$ be the smallest closed
$G_0$-invariant subalgebra containing $\{f_n:n\in \N\}$ and the
constant functions. Clearly $A$ is a separable subalgebra and the fact
that $G_0$ is dense in $G$ implies that $A$ is in fact $G$-invariant.

Let $Y$ be the compact metric Gelfand space of $A$.  (Thus the
elements of $Y$ are the multiplicative linear functionals of norm one
on $A$ and the map $A \cong C(Y),\ f \mapsto \hat f$, where $\hat
f(y)=y(f)$, is an isometric isomorphism of Banach algebras.)  Then,
for each $g\in G$, the linear action $ f \mapsto f\circ g $ of $g$ on
$A$ defines a homeomorphism $g: Y \to Y$ and $\widehat {f\circ g}
=\hat f\circ g$.  If $y_n\to y$ in $Y$ and $g_n\to e$ in $G$ are
convergent sequences then for every $\hat f\in C(Y)$ we have
\begin{align*}
|\hat{f}(g_n y_n)-\hat f(y)| & \le |\hat{f}(g_n y_n)-\hat f(y_n)| +
 |\hat f(y_n)-\hat f(y)| \\
& \le \|\widehat{f\circ g_n}-\hat f\| + |\hat f(y_n)-\hat f(y)|\,,
\end{align*}
hence $\lim_{n\to\infty}|\hat f(g_ny_n)-\hat f(y)|=0$.  It follows
that $\lim_{n\to\infty}g_ny_n=y$ and we conclude that the action of
$G$ on $Y$ is topological.

The linear functional $\mu: A \to \R, \ f\mapsto \int f\,d\mu$ defines
a probability measure $\nu$ on $Y$ and the dynamical system
$(Y,\Ycal,\nu,G)$, where $\Ycal$ is the Borel $\sig$-algebra on $Y$,
yields a Boolean action $(\Ycal,\nu,G)$ which is isomorphic to the
given Boolean action.  We conclude that $(Y,\Ycal,\nu,G)$ is a spatial
model as required.
\end{proof}

\begin{remark}\label{2.3}
In general when we do not assume that the $G$-continuous functions on
the near action $(X,\Xcal,\mu,G)$ separate points, we can still
consider the smallest $\sig$-algebra $\Dcal\subset \Xcal$ with respect
to which all the functions in $\Acal(G)$ are measurable and then the
closed subspace of $L^2(\mu)$ consisting of $\Dcal$-measurable
functions.  This subspace defines a factor near action and it is clear
that this factor is the largest factor which admits a spatial model.
\end{remark}

Theorem \ref{thm-levy} together with remark \ref{2.3} yield the
following:

\begin{cor}\label{levy}
For a L\'evy group $G$, an ergodic near-action admits only
constants as $G$-continuous functions.
\end{cor}

It is an interesting fact that a seemingly weaker condition already
implies $G$-continuity.  To see this we first need a lemma.

\begin{lem}\label{boco}
Let $ X $ be a Polish space, $ f : X \to L^1(\mu) $ a continuous map
such that the image $ f(X) $ is contained in $ L^\infty(\mu) $ and
that as a subset of the Banach space $ L^\infty(\mu) $ it is
separable.  Then $ f $ treated as a map $ X \to L^\infty(\mu) $ is
continuous at every point of some dense $ G_\delta $ subset of $ X $.
\end{lem}

\begin{proof}
Every closed ball in $ L^\infty(\mu) $ is a closed subset of $
L^1(\mu)$. We choose $ x_1, x_2, \dots \in X $ such that $ f(x_k) $
are $ L^\infty $-dense in $ f(X) $. We consider closed balls $ B_{n,r}
$ in $ L^\infty(\mu) $ of radius $ r $ centered at $ f(x_n) $. Their
inverse images $ f^{-1} (B_{n,r}) $ are closed in $ X $, and $ \cup_n
f^{-1} (B_{n,r}) = X $ (for every $ r>0 $). Denoting by $ U_{n,r} $
the interior of $ f^{-1} (B_{n,r}) $ we observe that $ \cup_n U_{n,r}
$ is a dense open set in $ X $ and $ \cap_r \cup_n U_{n,r} $ is a
dense $G_\delta $ set (Baire's theorem).

If $ y \in U_{n,r} $ and $ y_k \to y $ then $ f(y) \in B_{n,r} $ and $
f(y_k) \in B_{n,r} $ for large $ k $, therefore $ \limsup_k \| f(y_k)
- f(y) \|_\infty \le 2r $. So, if $ y \in \cap_r \cup_n U_{n,r} $ and
$ y_k \to y $ then $ \limsup_k \| f(y_k) - f(y) \|_\infty = 0 $.
\end{proof}

\begin{prop}\label{separable}
A function $ f \in L^\infty(\mu) $ is $G$-continuous if and only if
its $G$-orbit is a separable subset of $ L^\infty(\mu) $.
\end{prop}

\begin{proof}
The necessity is easy to see. The sufficiency follows from Lemma
\ref{boco}, applied to the map $ G \to L^\infty(\mu) $, $ g \mapsto f
\circ g $.  By homogeneity, its continuity at a single point implies
continuity everywhere.
\end{proof}

\section{Whirly actions}

Often one can use the necessary and sufficient condition of Theorem
\ref{2.2} to verify directly that a given near action has no spatial
model.  This is done most easily by the following notion which will
guarantee that a near action admits only constants as $G$-continuous
functions.

\begin{defn}\label{3.1}
A near-action of $G$ on $(X,\Xcal,\mu)$ is \emph{whirly,} if for all
sets $A,B\in\Xcal$ of positive measure, for almost all $g$ in $G$ with
respect to Baire category, $ \mu (A\cap gB) > 0 $.
\end{defn}

\begin{defn}\label{3.2}(equivalent to \ref{3.1}).
A near-action of $G$ on $(X,\Xcal,\mu)$ is \emph{whirly,} if for all
sets $A,B\in\Xcal$ of positive measure, every neighborhood of $e$ (the
unit of $G$) contains $g$ such that $ \mu (A\cap gB) > 0 $.
\end{defn}

Clearly, \ref{3.2} follows from \ref{3.1} (since a neighborhood cannot
be Baire-negligible). On the other hand, \ref{3.1} follows from
\ref{3.2}, since $ \mu(A\cap gB) $ is a continuous function of $ g $,
therefore the set $ V(A,B) = \{ g : \mu(A\cap gB) > 0 \} $ is
open. Its closure contains $e$. The same holds for the set $ V(gA,B) =
gV(A,B) $, which shows that $ V(A,B) $ is dense in $G$.

\begin{prop}\label{3.3}
\mbox{ }

\textup{(a)}
If a near-action is whirly, then all $G$-continuous functions
are constants.

\textup{(b)}
A whirly action has no spatial model; moreover, such an action cannot
have nontrivial spatial factors.
\end{prop}

\begin{proof}
(a)\
Assume that a $G$-continuous function $ f \in L^\infty(\mu) $ is
non-constant, then the sets $ A = f^{-1} ((-\infty,a)) $ and $ B =
f^{-1} ((b,+\infty)) $ are of positive measure, provided that $ a<b $
are chosen appropriately. All sufficiently small $g\in G$ (that is,
close enough to $e$) satisfy $ \| f - f \circ g^{-1} \|_\infty < b-a
$, therefore $ \mu(A\cap gB) = 0 $ and the action cannot be whirly.

(b)\
This follows from part 1 and Theorem \ref{2.2}.  The claim about the
factors follows from Remark \ref{2.3}.
\end{proof}

\begin{remark}\label{3.4}
Here is yet another equivalent definition.  A near-action of $G$ on
$(X,\Xcal,\mu)$ is whirly, iff for every set $A\in\Xcal$ of positive
measure and every neighborhood $U$ of $e$ in $G$,
\[
\mu (UA) = 1 \, ;
\]
here $ UA $ means $ \cup_n (g_n A) $ where $ (g_n) $ is a dense
sequence in $ U $ (its choice does not matter $ \bmod\,0 $).
Proof: $ \mu ((UA)\cap B)>0 \iff \exists n \; \mu((g_n A)\cap B)>0
\iff \exists g\in U \; \mu((gA)\cap B)>0 \iff \exists g\in U \;
\mu(A\cap g^{-1}B)>0 $.
\end{remark}

We will next describe some applications of Proposition \ref{3.3}.  Our
first application will be to the natural near-action (on $X$) of the
group $G=\Aut(X)$ of the automorphisms of the Lebesgue space
$(X,\Xcal,\mu)$.  We have already seen that this action has no spatial
model since $G$ is a L\'evy group. There is however a more direct
proof; we simply verify that the action is whirly.  To this end recall
that a neighborhood of the identity in $G$ is given by a finite
measurable partition of $X$ into sets $\{P_1,P_2,\dots,P_N\}$ and
$\ep>0$ as:
\[
U=\{S\in G: \sum_{j=1}^N \mu(P_j\tri SP_j)<\ep\}.
\]
For any sets $A,B \in \Xcal$ of positive measure, if $A_0\subset A,
B_0\subset B$ are measurable, disjoint and have the same measure
$\mu(A_0)=\mu(B_0)< \ep/2$, and $S$ is defined to be a measure
preserving transformation which is the identity on $X\setminus
(A_0\cup B_0)$ and interchanges $A_0$ with $B_0$, then $S\in U$ and it
satisfies $\mu(A \cap SB) > 0$.

The same kind of argument can be given for many subgroups of $G$ and
their natural near-action on $X$. For example we can start with any
countable subgroup $\Gamma\subset G$ that acts ergodically on $X$. The
\emph{full group} of this action $[\Gamma]$, consists of all the
measure preserving transformations $T\in G=\Aut(X,\Xcal,\mu)$ such
that for $\mu$ a.e. $x\in X$, $Tx\in \Gamma x$.

\begin{prop}\label{whirly}
The near action of $[\Ga]$ on $X$ is whirly.
\end{prop}

\begin{proof}
The argument given above for the entire group $G$ works here as
well, almost verbatim. The only place where some change is needed
is when we choose the transformation $S$; this time it should be
in $[\Ga]$.
Now suppose we are given the open set $U= U(P_0,P_1,\dots,P_N;\ep)$
and two positive measure sets $A,B \in X$. If $\mu(A\cap B)>0$ there
is nothing left to show. Otherwise the ergodicity of the $\Ga$-action
guarantees the existence of a $\ga\in \Ga$ with $\mu(A \cap \ga
B)>0$. Set $A_1=A \cap \ga B$ and choose any $A_0\subset A_1$ with $0
< \mu(A_0) <\ep/2$. We let $B_0=\ga^{-1} A_0$.  The transformation $S$
is now defined as the identity on $X\setminus (A_0\cup B_0)$ and it
interchanges $A_0$ with $B_0$ by means of $\ga$ and
$\ga^{-1}$. Clearly $S\in [\Ga]$, $S\in U$ and it satisfies $\mu(A
\cap SB) > 0$.
\end{proof}

\begin{rmk}
Let us note that, clearly, for every dense subgroup $H$ of
$G=\Aut(X,\Xcal,\mu)$ the action of $H$ on $X$ is whirly.  Moreover it
can be shown that in the notation of the previous discussion, the
group $[\Ga]$ is dense in $G$ whenever $\Ga$ acts ergodically on
$X$. However, proving the latter assertion requires a considerably more
elaborate argument than the direct proof we provided in Proposition
\ref{whirly}.  In addition to the topology of convergence in measure
on $G=\Aut(X,\Xcal,\mu)$ one can consider the much stronger
\emph{uniform topology} given by the metric $ d_u(S,T)=\mu\{x\in X:
Sx\ne Tx\}$.
Although with respect to this topology $G$ is a non Polish
topological group, the subgroup $[\Ga]$ is a closed Polish subgroup
(see Hamachi-Osikawa \cite{HO}, Lemma 53,
or observe directly that for
a countable generating collection of measurable partitions
$\Pcal=\{P=(P_1,\dots,P_n)\}$ of $X$ and a fixed enumeration
$\Ga=\{\ga_1,\ga_2,\dots\}$, the set of elements
$$
\{h(P;k_1,k_2,\dots,k_n) \in [\Ga]: P\in \Pcal,\ (k_1,k_2,\dots,k_n)
\in \N^n\},
$$
where
$$
h(P;k_1,k_2,\dots,k_n)\rest P_j =\ga_{k_j}, \ j=1,2,\dots,n,
$$
is a countable dense subset of $[\Ga]$).
We now note that for both  $G$ and $[\Ga]$ our proofs
show that their actions on $(X,\Xcal,\mu)$ are in fact whirly
with respect to the uniform topology.
Again we are indebted to V. Pestov for pointing this out.
\end{rmk}

\section{The automorphism group of the Gaussian measure}

We now turn back to the full orthogonal group $G=O(l^2) $ acting on
$(\R^\infty,\ga^\infty) $ as explained before Corollary
\ref{nogo1}. Thus we let $ \zeta_1, \zeta_2, \dots : X \to \R $ be
i.i.d.\ $ N(0,1) $ random variables defined as the coordinate
functions on the space of sequences $X=\R^\infty$ equipped with its
Borel $\sig$-algebra $\Xcal$ and the Gauss measure
$\mu=\ga^\infty$. We identify $ l^2 $ with the closed linear subspace
$H\subset L^2(\mu)$ generated by the functions
$\zeta_1,\zeta_2,\dots$; namely, $ (c_1,c_2,\dots) \in l^2 $ with
$c_1 \zeta_1 + c_2 \zeta_2 + \dots \in H $. The near-action is given by
\[
(c'_1 \zeta_1 + c'_2 \zeta_2 + \dots ) = ( c_1 \zeta_1 + c_2 \zeta_2 +
\dots ) \circ g \quad \text{whenever} \quad (c'_1,c'_2,\dots) =
g(c_1,c_2,\dots)
\]
for $ (c_1,c_2,\dots) \in l^2 $, $ g \in O(l^2) $; we call it the
automorphism group of the Gaussian measure.

\begin{thm}\label{Hample}
The near-action of the automorphism group of the Gaussian measure is
whirly.
\end{thm}

The following technical lemma is important for the proof. Basically it
states that a remote perturbation of a finite-dimensional condition
forces the conditional probability to be strictly positive. Note that
the relation $ \cP A \dots > 0 $ may be written as
$ \sgn \cP A \dots = 1 $, using the (discontinuous) sign function.

\begin{lem}\label{alpha}
Let $ A \subset X $ be a measurable set of positive probability, and $
\al \in (0,\pi/2) $. Then

\textup{(a)}
$ \sgn \cP{A}{ \zeta_1 \cos\al + \zeta_n \sin\al } \to 1 $ in
probability, for $n\to\infty $;

\textup{(b)}
for each $ m = 1,2,\dots $
\[
\sgn \cP{A}{ \zeta_1 \cos\al + \zeta_n \sin\al, \dots, \zeta_m \cos\al
 + \zeta_{n+m-1} \sin\al } \to 1 \quad \text{for } n \to \infty
\]
in probability.
\end{lem}

\begin{proof}
(a)\
We introduce functions $ f_n : \R \to [0,1] $,
$ g_n : \R^2 \to [0,1]$ by
\begin{gather*}
\cP{A}{ \zeta_1 \cos\al + \zeta_n \sin\al } = f_n ( \zeta_1 \cos\al +
 \zeta_n \sin\al ) \, , \\
\cP{A}{ \zeta_1, \zeta_n } = g_n ( \zeta_1, \zeta_n ) \, ,
\end{gather*}
and a set $ B \subset \R $ by
\[
B = \{ x_1 : \cP{A}{\zeta_1=x_1} > 0 \} \, .
\]
(These $ f_n, g_n, B $ are treated $ \bmod\,0 $, of course.) We have
\[
\cP{A}{\zeta_1,\zeta_n} \to \cP{A}{\zeta_1}
\]
in probability (for $ n \to \infty $). On $ B \times \R $ we get
\[
{\sgn} g_n \to 1 \quad \text{in measure,}
\]
with respect to  $ \ga \times \ga $, where $ \ga = N(0,1) $ is the
one-dimensional Gaussian measure. However, any equivalent (that is,
mutually absolutely continuous) finite measure on $ B \times \R $ may
be used equally well.

Taking into account that $ f_n $ results from $ g_n $ by integration
(along straight lines orthogonal to the unit vector $(\cos\al,
\sin\al)$) we get
\[
{\sgn} f_n (u) \ge {\ess\sup}_x \,{\sgn} g_n \bigg( x,
\frac{u-x\cos\al}{\sin\al} \bigg) \, .
\]
The map $ (x,u) \mapsto \( x, \frac{u-x\cos\al}{\sin\al} \) $ of $ B
\times \R $ to itself sends the measure to an equivalent measure. So,
$ {\sgn} g_n \( x, \frac{u-x\cos\al}{\sin\al} \) \to 1 $ in measure,
which implies $ {\sgn} f_n \to 1 $ in measure (with respect to $ \ga
$).

(b)\
The same as before, but $ \R $ is replaced by $ \R^m $, $ \R^2 $ by $
\R^{2m} $, $ \zeta_1 $ by $ (\zeta_1,\dots,\zeta_m) $ and $ \zeta_n $
by $ (\zeta_n,\dots,\zeta_{n+m-1}) $.
\end{proof}

\begin{proof}[Proof of Theorem \textup{\ref{Hample}}]
Let $ A \in \Xcal $ be a set of positive measure, and $ U $ a
neighborhood of $ e $ in $ G $; by Remark \ref{3.4} it is sufficient
to prove that $ \mu(UA) = 1 $. Of course, $ UA $ is treated as in
Remark \ref{3.4} (and the same about $ ZA $ for any $ Z \subset G $).

Ergodicity of $ G $ ensures that $ \mu (GA) = 1 $. Applying the same
argument to conditional measures we get (almost everywhere)
\[
\cP{ G_m A }{ \zeta_1,\dots,\zeta_m } \ge \sgn \cP{ A }{
\zeta_1,\dots,\zeta_m } \, ,
\]
where $ G_m = \{ g \in G : g \zeta_1 = \zeta_1, \dots, g \zeta_m =
\zeta_m \} $. For $ m $ large enough we have $ G_m \subset U $,
therefore
\[
\cP{ U A }{ \zeta_1,\dots,\zeta_m } \ge \sgn \cP{ A }{
\zeta_1,\dots,\zeta_m } \, .
\]
However, there is nothing special in $ \zeta_1,\dots,\zeta_m $; by the
$ O(l^2) $-invariance, the same holds for $ \xi_1 = \zeta_1 \cos\al +
\zeta_n \sin\al, \dots, \xi_m = \zeta_m \cos\al + \zeta_{n+m-1}
\sin\al $ provided that $ n>m $ and the corresponding subgroup $
G_{m,n,\al} = \{ g \in G : g \xi_1 = \xi_1, \dots, g \xi_m = \xi_m \}
$ is contained in $ U $. We choose $ m $ so large and $ \al $ so small
that $ G_{m,n,\al} \subset U $ for every $ n > m $ (this is possible
since for every $ h \in H $ its distance from the span of $ \xi_1,
\dots, \xi_m $ tends to $ 0 $ uniformly in $ n $ for $ m \to \infty $,
$ \al \to 0 $). We have
\[
\cP{ U A }{ \xi_1,\dots,\xi_m } \ge \sgn \cP{ A }{ \xi_1,\dots,\xi_m }
\]
for all $ n>m $. For $ n\to\infty $ the right-hand side converges to $
1 $ in probability (therefore, in $ L^1 $) by Lemma
\ref{alpha}. Taking the expectation we get $ \mu(UA) = 1 $.
\end{proof}

\appendix

\section{}

Measure concentration (that is, the property of being a L\'evy family
or group) is proven for various cases by a number of methods
\cite{Mi4}, \cite{L}.
Strong results need complicated proofs involving advanced
methods (Riemann geometry, representation theory, etc.). More
elementary arguments give weaker results which are satisfactory for
many topological applications such as the ones we needed in Section 1.
This appendix collects elementary (complete) proofs for many L\'evy
families.

\subsection{}\label{(1)}

Consider Gaussian measures $ \ga_\sig^n $ on $ \R^n $,
\[
\ga_\sig^n (dx) = (2\pi)^{-n/2} \sig^{-n} \exp \bigg( \! - \frac{
|x|^2 }{ 2\sig^2 } \bigg) \, dx \, ;
\]
note that $ \ga_\sig^n (\R^n) = 1 $ and $ \int |x|^2 \, \ga_\sig^n(dx)
= n \sig^2 $. We claim that $ ( \R^n, d_n, \ga_{\sig_n}^n ) $ is a
L\'evy family whenever the positive numbers $ \sig_n $ satisfy $ \sig_n
\to 0 $; here $ d_n (x,y) = |x-y| $ is the usual Euclidean metric on $
\R^n $. (Only the case $ \sig_n = n^{-1/2} $ will be used.)

According to the well-know relation between L\'evy families and
Lipschitz functions \cite[Sect.~1.3]{L}, it suffices to prove the
inequality
\[
\int f^2 \, d\ga_\sig^n \le \sig^2 \| f \|^2_\Lip
\]
for all functions $ f : \R^n \to \R $ such that $ \int f \,
d\ga_\sig^n = 0 $ and
\[
\| f \|_\Lip = \sup_{x\ne y} \frac{ |f(x)-f(y)| }{ |x-y| } < \infty \,
.
\]
Here is a proof. We introduce functions $ \phi, u : \R^n \times
(0,\sig^2) \to \R $,
\begin{align*}
\phi (x,t) = (2\pi)^{-n/2} t^{-n/2} \exp \bigg( \! - \frac{ |x|^2
 }{ 2t } \bigg) \, , \\
u (x,t) = \int f(y) \phi (y-x, \sig^2-t) \, dy \, ;
\end{align*}
they satisfy the (famous) partial differential equations
\[
\bigg( \frac{\partial}{\partial t} - \frac12 \Del \bigg) \phi(x,t) = 0
\, , \quad
\bigg( \frac{\partial}{\partial t} + \frac12 \Del \bigg) u(x,t) = 0 \,
;
\]
here $ \Del = \frac{\partial^2}{\partial x_1^2} + \dots +
\frac{\partial^2}{\partial x_n^2} $. Note that $ u(0,0+) = \int f \,
d\ga^n_\sig = 0 $ and $ u(x,\sig^2-) = f(x) $. It remains to prove
the inequality
\[
\int u^2 (x,t) \phi (x,t) \, dx \le t
\]
for $ 0 < t < \sig^2 $, assuming $ \| f \|_\Lip \le 1 $, that is, $ |
\nabla u(x,t) | \le 1 $ for all $ x, t $; here $ \nabla $ is the
gradient (in $ x $). For $ t \to 0+ $ the integral tends to $ u^2(0,0)
= 0 $. We have
\[
\frac{d}{dt} \int u^2 (x,t) \phi (x,t) \, dx =
\int \Big( -u \Del u + \frac12 \Del u^2 \Big) \phi \, dx = \int |
 \nabla u |^2 \phi \, dx \le 1 \, ,
\]
which completes the proof.

See also \cite[pp.~42, 49]{L}.

\subsection{}\label{(2)}

Euclidean spheres $ S^{n-1} = \{ x \in \R^n : |x| = 1 \} $ are a
L\'evy family, since a random point of $ S^{n-1} $ can be obtained
from a Gaussian random vector $ \xi \in \R^n $ distributed $
\ga^n_{1/\sqrt n} $ by the normalization map $ \xi \mapsto
\frac{\xi}{|\xi|} $; the map belongs to $ \Lip (2) $ as far as $ |\xi|
\ge 1 / 2 $. The other case, $ |\xi| < 1 / 2 $, may be ignored, since
its probability tends to $ 0 $ for $ n \to \infty $. The argument
works also when the radius $ r_n $
of the sphere is not just $ 1 $ but
satisfies $ r_n = o(\sqrt n \,) $. See also \cite[2.1, 2.3]{GM},
\cite[Prop.\ 2.10]{L},
\cite[$3\frac12.24$]{Gr}.

An alternative, comparably elementary way to \ref{(1)} and \ref{(2)}
is, first proving \ref{(2)} via the spectral gap of the Laplace
operator on the sphere \cite[4.2(a)]{GM}, \cite[Th.~3.1 and p.~49]{L}
and then deriving \ref{(1)} from \ref{(2)} \cite[p.~28]{L}.

\subsection{}\label{(3)}

The Stiefel manifolds $ W_2^n = \{ (x_1,x_2) \in S^{n-1} \times
S^{n-1} : \ip{x_1}{x_2} = 0 \} $ are a L\'evy family, since a random
point of $ W_2^n $ can be obtained from a $ 2n $-dimensional Gaussian
random vector $ (\xi_1,\xi_2) \in \R^n \oplus \R^n $ distributed $
\ga^n_{1/\sqrt n} \otimes \ga^n_{1/\sqrt n} $ by normalization,
subsequent orthogonalization $ (\xi_1,\xi_2) \mapsto (\xi_1, \xi_2 -
\ip{\xi_2}{\xi_1} \xi_1 ) $ and normalization again. The Lipschitz
property is ensured as far as $ | \ip{\xi_1}{\xi_2} | \le \ep_2
|\xi_1| |\xi_2| $ and $ |\xi_1|, |\xi_2| \in [1-\ep_2, 1+\ep_2] $,
where $ \ep_2 $ is an appropriate absolute constant. The other case
may be ignored, since its probability tends to $ 0 $ for $ n \to
\infty $. The orthonormalization commutes with the natural action of $
O(n) $; thus, $ O(n) $-invariance of the Gaussian measure ensures $
O(n) $-invariance of the measure on $ W_2^n $. The same argument (with
$ \ep_k $ in place of $ \ep_2 $) works for $ W_k^n = \{
(x_1,\dots,x_k) \in (S^{n-1})^k : \ip{x_i}{x_j} = 0 \text{ \small for
} 1\le i < j \le k \} $.
See \cite{Mi2}, \cite{Mi3}.
The proof in \cite[3.3]{GM} is somewhat less
elementary.

\subsection{}\label{(4)}

The full orthogonal group $ O(H) $ of a separable infinite-dimensional
Hilbert space $ H $ over $ \R $ is a L\'evy group.

\emph{Proof.} We
equip $ O(H) $ with a left-invariant metric and consider the subgroups $
O(n) = \{ g \in O(H) : g e_{n+1} = e_{n+1},  g e_{n+2} = e_{n+2},
\dots \} $ where $ e_1, e_2, \dots $ are a chosen orthonormal basis of
$ H $. Let $ A_n \subset O(n) $, $ \liminf m_n (A_n) > 0 $. Given $
\ep > 0 $, we take $ k $ and $ \del $ such that
\[
\| g e_1 - g' e_1 \|^2 + \dots + \| g e_k - g' e_k \|^2 < \del^2
\imp d(g,g') < \ep
\]
for all $ g, g' \in O(H) $. The maps $ T_n : O(n) \to W_k^n $, $ T_n(g) =
(g e_1, \dots, g e_k) $ satisfy
\[
\dist ( T_n(g), T_n(g') ) < \del \imp d(g,g') < \ep
\]
for all $ g,g' \in O(H) $. Therefore
\[
B_\ep (A_n) \supset T_n^{-1} ( B_\del (T_n(A_n)) ) \, .
\]
It remains to apply \ref{(3)}.

\subsection{}\label{(5)}

The commutative Polish group $ L_0 ([0,1],S^1) $ is a L\'evy group. It
consists of all equivalence classes of measurable functions $ [0,1]
\to S^1 $, where $ S^1 = \{ z\in\C : |z| = 1 \} $, and is in fact
monothetic \cite{G}. The following proof of its L\'evy property is
basically an extract from \cite[Sect.\ 1.6, 4.1]{L}. See also
\cite[p.~31]{L} and \cite{G}.

Let $ G $ be a commutative Polish group with a compatible invariant
metric $ d $, and $ \mu $ a Borel probability measure on $ G $; we
define
\[
\VarLip (\mu) = \sup \bigg\{ \sqrt{ \textstyle\int f^2 \, d\mu } : \|
f \|_\Lip \le 1 , \, \textstyle\int f \, d\mu = 0 \bigg\} \, .
\]
Clearly,
\[
\VarLip (\mu) \le \diam \supp (\mu)
\]
(the diameter of the support). It is easy to see that
\[
\VarLip (\mu) = \sup_f \frac{ \| f^2*\mu - (f*\mu)^2 \|_{\sup}^{1/2} }{
\| f \|_\Lip } \, ,
\]
where $ (f*\mu) (x) = \int f(x-y) \, \mu(dy) $, $ \| f \|_{\sup} =
\sup_{x\in G} |f(x)| $, and the squares are taken pointwise. For any
two measures $ \mu, \nu $
\begin{multline*}
\| f^2 * \mu * \nu - (f*\mu*\nu)^2 \|_{\sup} \le \\
\le \| ( f^2*\mu - (f*\mu)^2 ) * \nu \|_{\sup} +
 \| (f*\mu)^2*\nu - (f*\mu*\nu)^2 \|_{\sup} \le \\
\le \| f^2 * \mu - (f*\mu)^2 \|_{\sup} + \VarLip^2 (\nu) \| f*\mu
 \|^2_\Lip \le \\
\le \VarLip^2 (\mu) \| f \|^2_\Lip + \VarLip^2 (\nu) \| f \|^2_\Lip
 \, ,
\end{multline*}
thus
\[
\VarLip (\mu*\nu) \le \sqrt{ \VarLip^2 (\mu) + \VarLip^2 (\nu) } \, .
\]

The argument is applied to $ G = L_0 ((0,1),S^1) $ as follows. We
choose the $ L_1 $-metric
\[
d(x,y) = \int_0^1 | x(t) - y(t) | \, dt \, .
\]
For each $ n $ we consider the $ n $-dimensional compact subgroup
\[
K_n = K_{n,1} + \dots + K_{n,n} \subset G
\]
where $ K_{n,m} $ is the one-dimensional group of functions constant on
$ (\frac{m-1}n, \frac m n) $ and equal to $ 1 $ on $ (0,1) \setminus
(\frac{m-1}n, \frac m n) $. The corresponding invariant measures are
related by
\[
m_{K_n} = m_{K_{n,1}} * \dots * m_{K_{n,n}} \, .
\]
However,
\[
\VarLip ( m_{K_{n,m}} ) \le {\diam} ( K_{n,m} ) \le \frac2n \, ,
\]
therefore $ \VarLip^2 (m_{K_n}) \le n \cdot ( \frac2n )^2 \to 0 $ for
$ n \to \infty $. It remains to use the relation between L\'evy
families and Lipschitz functions mentioned in \ref{(1)}.

\subsection{}\label{(6)}

The Polish group $ G = \Aut ([0,1]) $ is a L\'evy group. It consists
of all equivalence classes of invertible transformations $ [0,1] \to
[0,1] $ preserving Lebesgue measure. Its L\'evy property may be proven
by the argument of \ref{(5)}, generalized to an arbitrary (not just
commutative) Polish group $ G $ with a compatible right-invariant
metric $ d $ (that is, $ d(g_1 h, g_2 h) = d(g_1,g_2) $). Still, $
\VarLip^2 (\mu*\nu) \le \VarLip^2 (\mu) + \VarLip^2 (\nu) $
where $ \VarLip (\mu) $ is \emph{defined} as $ \sup_f \| f^2*\mu -
(f*\mu)^2 \|_{\sup}^{1/2} / \| f \|_\Lip $, $ f*\mu $ is defined by $
(f*\mu) (x) = \int f(xy^{-1}) \, \mu(dy) $, and $ \mu*\nu $ is defined
by $ \int f \, d(\mu*\nu) = \int f(xy) \, \mu(dx) \nu(dy) $. However,
the inequality $ \VarLip (\mu) \le \diam \supp (\mu) $ need not hold,
since the map $ y \mapsto xy^{-1} $ need not be isometric. If the
metric $ d $ is bi-invariant (that is, $ d(g_1 h, g_2 h) = d(g_1,g_2)
= d(h g_1, h g_2) $), then $ \VarLip (\mu) \le \diam \supp (\mu) $.

We apply the argument to the group $ S_n $ of all permutations of $
\{1,\dots,n\} $ equipped with the Hamming metric
\[
d(g,h) = \frac{ \# \{ k : g(k) \ne h(k) \} }{ n } \, .
\]
Its invariant measure $ m_{S_n} $ is the convolution of $ n $
measures, each concentrated on transpositions (that is, $ g $ such
that $ d(g,e) \le 2/n $). Indeed, $ S_{n-1} $ is naturally embedded
into $ S_n $, and $ m_{S_n} = m_{S_{n-1}} * \mu $ where $ \mu $ is
distributed uniformly on transpositions of $ n $ and $ k $ for $ k =
1,\dots,n $. So, $ \VarLip^2 (m_{S_n}) \le n \cdot ( \frac4n )^2 \to 0
$ for $ n \to \infty $.

It remains to note that there exists a natural embedding of the
inductive limit group $\mathbf{S}=\lim_{n\to\infty} S_{2^n}$ as a
dense subgroup of $ G = \Aut ([0,1]) $.  Here the group $S_{2^n}$ is
embedded into $S_{2^{n+1}}$ as the subgroup of permutations $\hat\sig$
of $\{0,1,\dots,2^{n+1}-1\}$ of the form $\hat\sig(2k)=2\sig(k)$ and
$\hat\sig(2k+1)=2\sig(k)+1,\ k=0,1,\dots, 2^n-1,\ \sig\in S_{2^n}$.
Then, for each $n$, the group $S_{2^n}$ is identified as the subgroup
of $ G = \Aut ([0,1]) $ which consists of the transformations
permuting the $2^n$ dyadic sub-intervals of $[0,1]$ by translations.
It can be easily seen that the restriction of the uniform metric
$d_u(S,T)=\mu\{x\in X: Sx\ne Tx\}$ on $G$ to $S_{2^n}$ is the Hamming
metric.  Thus with respect to this metric and using the estimation for
$ \VarLip^2 (m_{S_n})$ we see that $\mathbf{S}$ is a L\'evy
group. Since the identity map from $(\mathbf{S},d_u)$ to $G$ is
continuous and since $\mathbf{S}$ is dense in $G$ we can finally
conclude that also $G$ is a L\'evy group.  See \cite{GP} and also
\cite[Corollary 4.3]{L}.

\section{}

The theorem stated below is well known and widely used.
However it seems that complete proofs are not
easily found. The proof we provide is
from Oxtoby \cite{O} where it is attributed to Ulam.

\begin{thm}
A Polish topological group $ G $ which admits a Borel
$\sig$-finite (either right or left) invariant measure
class is locally compact.
\end{thm}

\begin{proof}
Let $\mu$ be a $\sig$-finite measure
on $G$ such that
$$
\mu(B)>0 \iff \mu(gB)>0,\qquad \forall\ {\text{measurable}}\ B,
\forall g\in G.
$$
Let $A$ be a measurable subset of $G$ with
$0 < \mu(A) < \infty$. Since every Borel measure
on a Polish space is regular there exists a compact set
$K\subset A$ with $0 < \mu(K) < \infty$.
Let $H < G$ be the subgroup
of $G$ which is generated by $L=K \cup K^{-1}$. Clearly
$H=\cup\{L^n: n\in \N\}$ is a $\sig$-compact group.
If $G/H$ is uncountable then there are uncountably many
distinct cosets of $H$ in $G$ and in particular uncountably many
pairwise disjoint translations of $K$. This however
contradicts the $\sig$-finiteness of $\mu$ and we
conclude that $G/H$ is countable. Now Baire's theorem implies
that ${\inte}L^n \ne\emptyset$ for some $n$
and we conclude that $G$ is
locally compact.
\end{proof}

\begin{rmk}
In \cite[theorem 7.1]{Ma0} Mackey proved a more general theorem.
He showed that if $G$ is an analytic Borel group
(no topology is given) which admits an invariant $\sig$-finite
measure class then there exists a unique locally compact topology
on $G$ whose Borel structure is the given one and under
which $G$ is a topological group. His proof relies on
Weil's theorem, \cite{We}.
Finally we note that in \cite{A} A. D. Alexandroff proved a
related result (he is not assuming $\sig$-finiteness of
the invariant measure but, instead, the existence of
an open set with finite positive measure).
\end{rmk}

\end{document}